\documentclass{amsart}

\usepackage{amsmath} 
\usepackage{amssymb}

\newtheorem{theorem}{Theorem}[section] 
\newtheorem{claim}{Claim}[theorem]
\newtheorem{lemma}[theorem]{Lemma} 
\newtheorem{proposition}[theorem]{Proposition} 
\newtheorem{observation}[theorem]{Observation} 
\newtheorem{corollary}[theorem]{Corollary} 

\theoremstyle{definition}
\newtheorem{definition}[theorem]{Definition}

\theoremstyle{remark}
\newtheorem{remark}[theorem]{Remark}

\numberwithin{equation}{section}
\newcommand{\lh}{{\ell g}}
\newcommand{\rest}{{\restriction}}
\newcommand{\dom}{{\rm dom}} 

\newcommand{\set}{{\rm set}}

\newcommand{\cB}{{\mathcal B}}

\newcommand{\gb}{{\mathfrak b}}

\newcommand{\cI}{{\mathcal I}}

\newcommand{\cM}{{\mathcal M}}
\newcommand{\gm}{{\mathfrak m}}

\newcommand{\bbP}{{\mathbb P}}

\newcommand{\gp}{{\mathfrak p}}

\newcommand{\bbQ}{{\mathbb Q}}

\newcommand{\bS}{{\bf S}}

\newcommand{\gt}{{\mathfrak t}} 
\newcommand{\cU}{{\mathcal U}}

\newcommand{\cf}{{\rm cf}}
\newcommand{\pr}{{\rm pr}}

\newcount\skewfactor
\def\mathunderaccent#1#2 {\let\theaccent#1\skewfactor#2
\mathpalette\putaccentunder}
\def\putaccentunder#1#2{\oalign{$#1#2$\crcr\hidewidth
\vbox to.2ex{\hbox{$#1\skew\skewfactor\theaccent{}$}\vss}\hidewidth}}
\def\name{\mathunderaccent\tilde-3 }

\begin{document}

\title{A comment on ``${\mathfrak p} < {\mathfrak t}$''}
\author{Saharon Shelah}
\address{Einstein Institute of Mathematics\\
Edmond J. Safra Campus, Givat Ram\\
The Hebrew University of Jerusalem\\
Jerusalem, 91904, Israel\\
 and  Department of Mathematics\\
 Rutgers University\\
 New Brunswick, NJ 08854, USA}
\email{shelah@math.huji.ac.il}
\urladdr{http://shelah.logic.at}
\thanks{The author acknowledges support from the United States-Israel
Binational Science Foundation (Grant no. 2002323). Publication 885.}     

\date{December 2006}

\begin{abstract}
Dealing with the cardinal invariants ${\mathfrak p}$ and ${\mathfrak t}$ of
the continuum we prove that 
\[{\mathfrak m}={\mathfrak p} = \aleph_2\ \Rightarrow\ {\mathfrak t} = 
\aleph_2.\]
In other words, if ${\bf MA}_{\aleph_1}$ (or a weak version of this) holds,
then (of course $\aleph_2\le {\mathfrak p}\le {\mathfrak t}$ and)
${\mathfrak p}=\aleph_2\ \Rightarrow\ {\mathfrak p}={\mathfrak t}$.  The
proof is based on a criterion for ${\mathfrak p}<{\mathfrak t}$. 
\end{abstract}

\maketitle

\section{Introduction}
We are interested in two cardinal invariants of the continuum, ${\mathfrak
p}$ and ${\mathfrak t}$. The cardinal ${\mathfrak p}$ measures when a family
of infinite subsets of $\omega$ with finite intersection property has a
pseudo--intersection.  A relative is ${\mathfrak t}$, which deals with 
towers, i.e., families well ordered by almost inclusion.  These are
closely related classical cardinal invariants. Rothberger \cite{Ro39},
\cite{Ro48} proved (stated in our terminology) that ${\mathfrak p}\leq
{\mathfrak t}$ and 
\[{\mathfrak p}=\aleph_1\ \Rightarrow\ {\mathfrak p} = {\mathfrak t},\]
and he asked if ${\mathfrak p} = {\mathfrak t}$.

Our main result is  Corollary \ref{y.49} stating that 
\[{\mathfrak m}\geq{\mathfrak p}=\aleph_2\ \Rightarrow\ {\mathfrak p}=
{\mathfrak t},\]
where ${\mathfrak m}$ is the minimal cardinal $\lambda$ such that Martin
Axiom for $\lambda$ dense sets fails (i.e. $\neg{\bf MA}_\lambda$). 
Considering that ${\mathfrak m}\ge\aleph_1$ is a theorem (of ZFC), the
parallelism with Rothberger's theorem is clear.  The reader may conclude
that probably ${\mathfrak m}={\mathfrak p}\ \Rightarrow\ {\mathfrak p} =
{\mathfrak t}$; this is not unreasonable but we believe that eventually one
should be able to show ${\rm CON}({\mathfrak m}=\lambda+{\mathfrak p}=
\lambda+{\mathfrak t}=\lambda^+)$. In the first section we present a
characterization of ${\mathfrak p}<{\mathfrak t}$ which is crucial for the
proof of \ref{y.49}, and which also sheds some light on the strategy to
approach the question of $\gp<\gt$ presented in \cite{Sh:F769}. 
\medskip

We thank Andreas Blass for detailed comments on an earlier version and
David Fremlin for historical information.
\medskip

\noindent{\bf Notation}\qquad Our notation is rather standard and compatible
with that of classical textbooks (like Bartoszy\'nski and Judach
\cite{BaJu95}). In forcing we keep the older convention that {\em a stronger
condition is the larger one}.  

\begin{enumerate}
\item Ordinal numbers will be denoted be the lower case initial letters of
the Greek alphabet ($\alpha,\beta,\gamma,\delta\ldots$) and also by $i,j$
(with possible sub- and superscripts). 
\item Cardinal numbers will be called $\kappa,\kappa_i,\lambda$. 
\item A bar above a letter denotes that the considered object is a sequence;
usually $\bar{X}$ will be $\langle X_i:i<\zeta\rangle$, where $\zeta$
is the length $\lh(\bar{X})$ of $\bar{X}$. Sometimes our sequences will be
indexed by a set of ordinals, say $S\subseteq\lambda$, and then $\bar{X}$
will typically be $\langle X_\delta:\delta\in S\rangle$. 
\item The set of all infinite subsets of the set $\omega$ of natural numbers
  is denoted by $[\omega]^{\aleph_0}$ and the relation of {\em almost
    inclusion\/} on $[\omega]^{\aleph_0}$ is denoted by $\subseteq^*$. Thus
  for $A,B\in [\omega]^{\aleph_0}$ we write $A\subseteq^* B$ if and only if
  $A\setminus B$ is finite. 
\item The relations of {\em eventual dominance\/} on the Baire space
  ${}^\omega \omega$ are called $\leq^*$ and $<^*$. Thus, for $f,g\in
  {}^\omega \omega$,
\begin{itemize}
\item $f\leq^*g$ if and only if $(\forall^\infty n<\omega)(f(a)\leq g(n))$ and   
\item $f<^*g$ if and only if $(\forall^\infty n<\omega)(f(a)< g(n))$.  
\end{itemize}
\end{enumerate}

\section{A criterion}
In this section our aim is to prove Theorem \ref{y.10} stating that
${\mathfrak p}<{\mathfrak t}$ implies the existence of a peculiar cut in
$({}^\omega \omega,<^*)$.  This also gives the background for our tries to
get a progress on the consistency of ${\mathfrak p}<{\mathfrak t}$ in
\cite{Sh:F769}.    

\begin{definition}
\begin{enumerate}
\item We say that a set $A\in [\omega]^{\aleph_0}$ is  {\em a
pseudo--intersection of a family $\cB\subseteq
[\omega]^{\aleph_0}$} if $A\subseteq^* B$ for all $B\in\cB$.  
\item A sequence $\langle X_\alpha:\alpha<\kappa\rangle\subseteq
  [\omega]^{\aleph_0}$ is {\em a tower\/} if $X_\beta\subseteq^* X_\alpha$
  for $\alpha<\beta<\kappa$ but the family $\{X_\alpha:\alpha<\kappa\}$ has
  no pseudo--intersection. 
\item $\gp$ is the minimal cardinality of a family $\cB\subseteq
[\omega]^{\aleph_0}$ such that the intersection of any finite subcollection
of $\cB$ is infinite but $\cB$ has no pseudo--intersection, and $\gt$ is the
smallest size of a tower. 
  \end{enumerate}
\end{definition}

A lot of results have been accumulated on these two cardinal invariants. For
instance 
\begin{itemize}
\item Bell \cite{Bell81} showed that $\gp$ is the first cardinal $\mu$  for
  which ${\bf MA}_\mu(\sigma\mbox{-centered})$ fails, 
\item Szyma\'nski proved that ${\mathfrak p}$ is regular (see, e.g., Fremlin
  \cite[Proposition 21K]{Fre}),
\item Piotrowski and Szyma\'nski \cite{PiSz87} showed that $\gt\leq {\rm
    add}({\cM})$ (so also $\gt\leq \gb$).
\end{itemize}
For more results and discussion we refer the reader to \cite[\S1.3,
\S2.2]{BaJu95}.  

\begin{definition}
\label{y.1} 
We say that a family $\cB\subseteq [\omega]^{\aleph_0}$ {\em exemplifies
  ${\gp}$} if:
\begin{itemize}
\item $\cB$ is closed under finite intersections (i.e., $A,B \in \cB\
  \Rightarrow\ A \cap B \in \cB$), and 
\item $\cB$ has no pseudo--intersection and $|\cB|=\gp$.
\end{itemize}
\end{definition}

\begin{proposition}
\label{y.2}
Assume $\gp<\gt$ and let ${\cB}$ exemplify $\gp$. Then there are a cardinal
$\kappa=\cf(\kappa) < \gp$ and a $\subseteq^*$--decreasing sequence $\langle
A_i:i < \kappa \rangle\subseteq [\omega]^{\aleph_0}$ such that
\begin{enumerate}
\item[(a)] $A_i \cap B$ is infinite for every $i<\kappa$ and $B\in\cB$, and  
\item[(b)] if $A$ is a pseudo--intersection of  $\{A_i:i<\kappa\}$, then for
  some $B\in\cB$ the intersection $A\cap B$ is finite. 
\end{enumerate}
\end{proposition}

\begin{proof}
Fix an enumeration $\cB=\{B_i:i<\gp\}$. By induction on $i < \gp$ we try to
choose $A_i \in [\omega]^{\aleph_0}$ such that 
\begin{enumerate}
\item[$(\alpha)$]  $A_i \subseteq^* A_j$ whenever $j<i$, 
\item[$(\beta)$]  $B \cap A_i$ is infinite for each $B\in\cB$, 
\item[$(\gamma)$] if $i=j+1$, then $A_i \subseteq B_j$.
\end{enumerate}
If we succeed, then $\{A_i:i < \gp\}$ has no pseudo--intersection so
${\mathfrak t} \leq\gp$, a contradiction.  So for some $i<\gp$ we cannot
choose $A_i$.  Easily $i$ is a limit ordinal and let $\kappa = \cf(i)$ (so
$\kappa \le i<\gp$). Pick an increasing sequence $\langle
j_\varepsilon:\varepsilon < \kappa \rangle$ with limit $i$.  Then $\langle
A_{j_\varepsilon}:\varepsilon < \kappa \rangle$ is as required.
\end{proof}

\begin{remark}
Concerning Proposition \ref{y.2}, let us note that Todor\v{c}evi\'{c} and
Veli\v{c}kovi\'{c} used this idea in \cite[Thm 1.5]{ToVe87} to exhibit a
$\sigma$--linked poset of size $\gp$ that is not $\sigma$--centered. 
\end{remark}

\begin{lemma}
\label{y.4}
Assume that
\begin{enumerate}
\item[(i)]  $\bar{A}=\langle A_i:i < \delta\rangle$ is a sequence of members
  of $[\omega]^{\aleph_0}$, $\delta<\gt$, 
\item[(ii)]  $\bar{B}=\langle B_n:n<\omega\rangle\subseteq
  [\omega]^{\aleph_0}$ is $\subseteq^*$--decreasing and 
\item[(iii)] for each $i<\delta$ and $n<\omega$ the intersection $A_i \cap
  B_n$ is infinite, and 
\item[(iv)] $\big(\forall i<j<\delta\big)\big(\exists n<\omega\big)\big(A_j
  \cap B_n \subseteq^* A_i\cap B_n\big)$. 
\end{enumerate}
Then for some $A \in [\omega]^{\aleph_0}$ we have 
\[\big(\forall i<\delta\big)\big(A \subseteq^* A_i\big)\mbox{ and }
\big(\forall n<\omega\big)\big(A\subseteq^* B_n\big).\]    
\end{lemma}

\begin{proof}
Without loss of generality $B_{n+1} \subseteq B_n$ and $\emptyset =
\bigcap\{B_n:n < \omega\}$ (as we may use $B'_n = \bigcap\limits_{\ell \le
  n} B_\ell \setminus \{0,\ldots,n\}$). For each $i < \delta$ let $f_i
\in {}^\omega \omega$ be defined by 
\[f_i(n)=\min\{k \in B_n \cap A_i: k > f_i(m)\mbox{ for every }m<n\}+1.\]
Since $\gt\leq\gb$, there is $f\in {}^\omega\omega$ such that $\big(\forall
i<\kappa\big)\big( f_i<^* f\big)$ and $n<f(n)< f(n+1)$ for $n<\omega$.
Let 
\[B^*=\bigcup\{(B_{n+1}\cap [n,f(n+1)):n < \omega\}.\]
Then $B^*\in [\omega]^{\aleph_0}$ as for $n$ large enough, 
\[\min[B_{n+1}\setminus [0,n) \cap A_0]\leq f_0(n+1)< f(n+1).\]
Clearly  for each $n<\omega$ we have  $B^* \setminus [0,f(n)) \subseteq B_n$
and hence $B^* \subseteq^* B_n$. Moreover, $(\forall i<\kappa)(A_i \cap B^*
\in [\omega]^{\aleph_0})$ (as above) and $(\forall i<j<\kappa)(A_j \cap B^*
\subseteq^* A_i\cap B^*)$ (remember assumption (iv)). Now applying $\gt
>\delta$ to $\langle A_i \cap B^*:i < \delta\rangle$ we get a
pseudo--intersection $A^*$ which is as required.
\end{proof}

\begin{definition}
\label{y.5}
\begin{enumerate}
\item Let $\bS$ be the family of all sequences $\bar{\eta}=\langle\eta_n:n
\in B\rangle$ such that $B\in [\omega]^{\aleph_0}$, and for $n\in B$,
$\eta_n \in {}^{[n,k)}2$ for some $k \in (n,\omega)$. We let
$\dom(\bar{\eta})=B$ and let $\set(\bar{\eta})=\bigcup\{\set(\eta_n):n \in
\dom(\bar{\eta})\}$,  where $\set(\eta_n)=\{\ell:\eta_n(\ell)=1\}$.
\item For $\bar{A}=\langle A_i:i<\alpha\rangle\subseteq[\omega]^{\aleph_0}$
  let 
\[\bS_{\bar{A}} =\big\{\bar{\eta}\in\bS:\big(\forall
i<\alpha\big)\big(\set(\bar{\eta})\subseteq^*A_i\big)\mbox{ and }
\big(\forall n\in\dom(\bar{\eta})\big)\big(\set(\eta_n) \ne \emptyset \big)
\big\}.\]    
\item For $\bar{\eta},\bar{\nu}\in\bS$ let $\bar{\eta}\leq^*\bar{\nu}$ mean
  that for every $n$ large enough, 
\[n\in\dom(\bar{\nu})\ \Rightarrow\ n \in\dom(\bar{\eta}) \wedge \eta_n
\trianglelefteq\nu_n\]
(where $\eta_n\trianglelefteq\nu_n$ means ``$\eta_n$ is an initial segment
of $\nu_n$'').   
\item For $\bar{\eta},\bar{\nu}\in\bS$ let $\bar{\eta}\leq^{**} \bar{\nu}$
  mean that for every $n \in \dom(\bar{\nu})$ large enough, for some $m \in
  \dom(\bar{\eta})$ we have $\eta_m\subseteq\nu_n$ (as functions).
\item For $\bar{\eta}\in\bS$ let  $C_{\bar{\eta}}=\{\nu\in {}^\omega
  2:(\exists^\infty n)(\eta_n\subseteq\nu)\}$.
\end{enumerate}
\end{definition}

\begin{observation}
\label{y.6A}
\begin{enumerate}
\item If $\bar{\eta} \le^* \bar{\nu}$, then $\bar{\eta}\le^{**} \bar{\nu}$
    which implies $C_{\bar{\nu}} \subseteq C_{\bar{\eta}}$. 
\item For every $\bar{\eta}\in\bS$ and a meagre set $B\subseteq {}^\omega 
2$ there is $\bar{\nu}\in\bS$ such that $\bar{\eta} \le^* \bar{\nu}$ and
$C_{\bar{\nu}}\cap B= \emptyset$.
\end{enumerate}
\end{observation}

\begin{lemma}
\label{y.6}
\begin{enumerate}
\item If $\bar{A}=\langle A_i:i<i^*\rangle\subseteq[\omega]^{\aleph_0}$
  has finite intersection property and $i^*<\gp$, then $\bS_{\bar{A}}
  \ne \emptyset$. 
\item Every $\le^*$-increasing sequence of members of $\bS$ of length $<\gt$ 
  has an $\le^*$-upper bound.    
\item If $\bar{A}=\langle A_i:i<i^*\rangle\subseteq[\omega]^{\aleph_0}$ is
  $\subseteq^*$--decreasing and $i^*<\gp$, then every $\le^*$-increasing
  sequence of members of $\bS_{\bar{A}}$ of length $<\gp$ has an
  $\le^*$-upper bound in $\bS_{\bar{A}}$.
\end{enumerate}
\end{lemma}

\begin{proof}
  (1)\quad Let $A \in [\omega]^{\aleph_0}$ be such that $(\forall i<i^*)(A
  \subseteq^* A_i)$ (exists as $i^*<\gp$).  Let $k_n=\min(A\setminus
  (n+1))$, and let $\eta_n \in {}^{[n,k_n+1)} 2$ be defined by
  $\eta_n(\ell)$ is $0$ if $\ell \in [n,k_n)$ and is $1$ if $\ell=k_n$. Then
  $\langle\eta_n: n<\omega\rangle\in \bS_{\bar{A}}$.  \smallskip

\noindent (2)\quad Let $\langle\bar{\eta}^\alpha:\alpha< \delta\rangle$ be a
$\le^*$--increasing sequence and $\delta<{\mathfrak t}$.  Let $A^*_\alpha =:
\dom(\bar{\eta}^\alpha)$ for $\alpha < \delta$. Then $\langle A^*_\alpha:
\alpha<\delta\rangle$ is a $\subseteq^*$--decreasing sequence of members of 
$[\omega]^{\aleph_0}$.  As $\delta < {\mathfrak t}$ there is $A^* \in
[\omega]^{\aleph_0}$ such that $\alpha < \delta\ \Rightarrow\ A^*\subseteq^*
A_\alpha^*$.  Now for $n < \omega$ we define 
\[B_n=\bigcup\big\{{}^{[m,k)}2: m \in A^*\mbox{ and }n \le m < k <
\omega\big\},\] 
and for $\alpha < \delta$ we define 
\[A_\alpha=\big\{\eta:\mbox{ for some }n\in\dom(\bar{\eta}^\alpha)\mbox{ we
  have }\eta^\alpha_n \trianglelefteq \eta\big\}.\]
One easily verifies that the assumptions of Proposition \ref{y.4} are
satisfied (well, replacing $\omega$ by $B_0$!!).  Let $A\subseteq B_0$ be
given by the conclusion of \ref{y.4}, and put  
\[A' =\big\{n:\mbox{ for some $\eta \in A$ we have }\eta \in
\bigcup\{{}^{[n,k)}2:k \in (n,\omega)\}\big\}.\]
Plainly, the set $A'$ is infinite. We let $\bar{\eta}^*=\langle\eta_n:n\in
A'\rangle$ where $\eta_n$ is any member of $A \cap B_n \setminus B_{n+1}$.   
\smallskip

\noindent(3)\quad Assume that $\bar{A}=\langle A_i:i<i^*\rangle \subseteq
[\omega]^{\aleph_0}$ is $\subseteq^*$--decreasing, $i^*<\gp$, and
$\langle\bar{\eta}^\alpha: \alpha< \delta\rangle\subseteq\bS_{\bar{A}}$ is
$\le^*$--increasing, and $\delta<\gp$. Let us consider the following forcing
notion $\bbP$.

\noindent {\bf A condition in $\bbP$} is a quadruple $p=(\bar{\nu},u,w,a)= 
(\bar{\nu}^p,u^p,w^p,a^p)$ such that 
\begin{enumerate}
\item[(a)] $u\in [\omega]^{<\aleph_0}$, $\bar{\nu}=\langle \nu_n:n\in
  u\rangle$, and for $n\in u$ we have:
\begin{itemize}
\item  $\nu_n\in {}^{[n,k_n)}2$ for some $k_n\in (n,\omega)$, and 
\item $\set(\nu_n)\neq \emptyset$,
\end{itemize}
\item[(b)] $w\subseteq \delta$ is finite, and  
\item[(c)] $a\subseteq i^*$ is finite.
\end{enumerate}

\noindent {\bf The order $\leq_{\bbP}=\leq$ of $\bbP$} is given by $p\leq q$
if and only if ($p,q\in\bbP$ and) 
\begin{enumerate}
\item[(i)]   $u^p\subseteq u^q$, $w^p\subseteq w^q$, $a^p\subseteq a^q$, and
  $\bar{\nu}^q\rest u^p=\bar{\nu}^p$, 
\item[(ii)]   $\max(u^p)<\min(u^q\setminus u^p)$ and for $n\in u^q\setminus
  u^p$ we have
\begin{enumerate}
\item[(iii)] $(\forall \alpha\in w^p)(n\in \dom(\bar{\eta}^\alpha)\ \wedge\
  \eta^\alpha_n\vartriangleleft \nu^q_n)$, 
\item[(iv)] $(\forall i\in a^p)(\set(\nu^q_n)\subseteq A_i)$. 
\end{enumerate}
\end{enumerate}

Plainly, $\bbP$ is a $\sigma$--centered forcing notion and the sets
\[\cI^{\alpha,i}_m=\big\{p\in\bbP: \alpha\in w^p\ \wedge \  i\in a^p\
\wedge\ |u^p|>m\big\}\]
(for $\alpha<\delta$, $i<i^*$ and $m<\omega$) are open dense in $\bbP$.
Since $|\delta|+|i^*|+\aleph_0<\gp$, we may choose a directed set $G
\subseteq\bbP$ meeting all the sets $\cI^{\alpha,i}_n$. Putting $\bar{\nu}=
\bigcup\{\bar{\nu}^p:p\in G\}$ we will get an upper bound to $\langle
\bar{\eta}^\alpha: \alpha< \delta\rangle$ in $\bS_{\bar{A}}$.  
\end{proof}

\begin{lemma}
\label{y.9}  
Assume that
\begin{enumerate}
\item[(i)] $\gp<\gt$ and ${\mathcal B}=\{B_\alpha:\alpha<\gp\}$ exemplifies
  $\gp$ (see \ref{y.1}), and  
\item[(ii)] $\bar{A}=\langle A_i:i<\kappa\rangle\subseteq
  [\omega]^{\aleph_0}$ is $\subseteq^*$--decreasing, $\kappa < \gp$ and
  conditions (a)+(b) of  Proposition \ref{y.2} hold, 
\item[(iii)]  $\pr:\gp\times\gp\longrightarrow\gp$ is a bijection
  satisfying  $\pr(\alpha_1,\alpha_2)\ge\alpha_1,\alpha_2$.   
\end{enumerate}
Then we can find a sequence $\langle\bar{\eta}^\alpha:\alpha\le\gp\rangle$
such that  
\begin{enumerate}
\item[(a)]  $\bar{\eta}^\alpha\in\bS_{\bar{A}}$ for $\alpha<\gp$ and
  $\bar{\eta}^{\gp}\in \bS$ (sic!), 
\item[(b)] $\langle\bar{\eta}^\alpha:\alpha\le\gp\rangle$ is
  $\le^*$-increasing,  
\item[(c)] if $\alpha<\gp$ and $n\in\dom(\bar{\eta}^{\alpha+1})$ is large
  enough, then $\set(\eta^{\alpha+1}_n)\cap B_\alpha\ne\emptyset$ (hence
  $(\forall^\infty n\in\dom(\bar{\eta}^\beta))(\set(\eta^\beta_n)\cap
  B_\alpha\ne\emptyset)$ holds for every $\beta\in [\alpha+1,\gp]$),  
\item[(d)] if $\alpha=\pr(\beta,\gamma)$, then $\set(\eta^{\alpha +1}_n)\cap
  B_\beta \ne \emptyset$ and $\set(\eta^{\alpha +1}_n) \cap B_\gamma\ne
  \emptyset$ for $n\in\dom(\bar{\eta}^{\alpha+1})$, and the truth values of
\[\min(\set(\eta^{\alpha+1}_n) \cap B_\beta) <\min(\set(\eta^{\alpha+1}_n)
  \cap B_\gamma)\]
are the same for all $n\in\dom(\bar{\eta}^{\alpha+1})$,
\item[(e)] in (d), if $\beta < \kappa$ we can replace $B_\beta$ by
  $A_\beta$; similarly with $\gamma$; and if $\beta,\gamma < \kappa$ then we
  can replace both. 
\end{enumerate}
\end{lemma}

\begin{proof}
We choose $\bar{\eta}^\alpha$ by induction on $\alpha$. For $\alpha = 0$
it is trivial, for $\alpha$ limit $< \gp$ we use Lemma \ref{y.6}(3) (and $|
\alpha|<\gp$). At a successor stage $\alpha+1$, we let $\beta,\gamma$ be
such that $\pr(\beta,\gamma)=\alpha$ and we choose $B'_\alpha \in [\omega
]^{\aleph_0}$ such that $B'_\alpha \subseteq B_\alpha\cap B_\beta\cap
B_\gamma$ and $(\forall i<\kappa)(B'_\alpha \subseteq^* A_i)$.  Next, for
$n\in \dom(\bar{\eta}^\alpha)$ we choose $\eta'_n$ such that
$\eta^\alpha_n\vartriangleleft \eta'_n$ and 
\[\emptyset\neq \{\ell:\eta'_n(\ell)=1\mbox{ and }\lh(\eta^\alpha_n)\le \ell
< \lh(\eta'_n)\}\subseteq B'_\alpha.\] 
Then we let $\bar{\eta}^{\alpha+1}=\langle\eta'_n:n\in\dom(
\bar{\eta}^\alpha)\rangle$.  By shrinking the domain of
$\bar{\eta}^{\alpha+1}$ there is no problem to take care of clause (d). It
should be also clear that me may ensure clause (e) as well.  

For $\alpha=\gp$, use \ref{y.6}(2).  
\end{proof}

\begin{definition}
\label{cut}
Let $\kappa_1,\kappa_2$ be infinite regular cardinals. {\em A
  $(\kappa_1,\kappa_2)$--peculiar cut in ${}^\omega\omega$} is a pair
  $\big(\langle f_i:i<\kappa_1\rangle,\langle
  f^\alpha:\alpha<\kappa_2\rangle\big)$ of sequences of functions in
  ${}^\omega \omega$ such 
\begin{enumerate}
\item[$(\alpha)$]   $(\forall i<j<\kappa_1)(f_j<^* f_i)$, 
\item[$(\beta)$]   $(\forall\alpha<\beta<\kappa_2)(f^\alpha<^* f^\beta)$,
\item[$(\gamma)$]   $(\forall i<\kappa_1)(\forall\alpha<\kappa_2)(f^\alpha
  <^* f_i)$, 
\item[$(\delta)$]  if $f:\omega\longrightarrow\omega$ is such that $(\forall
  i<\kappa_1)(f \le^* f_i)$, then $f \le^* f^\alpha$ for some
  $\alpha<\kappa_2$,  
\item[$(\varepsilon)$]  if $f:\omega\longrightarrow\omega$ is such that
  $(\forall \alpha<\kappa_2)(f^\alpha \le^* f)$, then $f_i \le^* f$ for some
  $i<\kappa_1$. 
\end{enumerate}
\end{definition}

\begin{proposition}
\label{getb}
If $\kappa_2<\gb$, then there is no $(\aleph_0,\kappa_2)$--peculiar cut.
\end{proposition}

\begin{proof}
Assume towards contradiction that $\gb>\kappa_2$ but there is an
$(\aleph_0,\kappa_2)$--peculiar cut, say $\big(\langle f_i:i<\omega\rangle,
\langle f^\alpha:\alpha<\kappa_2\rangle\big)$ is such a cut. Let $S$ be the
family of all increasing sequences $\bar{n}=\langle n_i:i<\omega\rangle$
with $n_0=0$. For $\bar{n}\in S$ and $g\in {}^\omega\omega$ we say that {\em
  $\bar{n}$ obeys $g$} if $(\forall i<\omega)(g(n_i)<n_{i+1})$. Also for
$\bar{n}\in S$ define $h_{\bar{n}}\in {}^\omega\omega$ by 
\[h_{\bar{n}}\restriction [n_i,n_{i+1})=f_i\restriction [n_i,n_{i+1})\qquad
    \mbox{ for }i<\omega.\] 
Now, let $g^*\in {}^\omega\omega$ be an increasing function such that for
every $n<\omega$ and $m\geq g^*(n)$ we have 
\[f_{n+1}(m)<f_n(m)<\ldots<f_1(m)<f_0(m).\]
Note that 
\begin{enumerate}
\item[$(\boxdot)_1$] if $\bar{n}\in S$ obeys $g^*$, then $(\forall
  i<\omega)(h_{\bar{n}}<^* f_i)$.
\end{enumerate}
Now, for $\alpha<\kappa_2$ define $g^\alpha\in {}^\omega\omega$ by  
\begin{enumerate}
\item[$(\boxdot)_2$] $g^\alpha(n)=\min\big\{k<\omega: k>n+1\ \wedge\ \big
  (\forall i\leq n\big)\big(\exists\ell\in [n,k)\big)\big(f^\alpha(\ell)<
  f_i(\ell)\big)\big\}$.
\end{enumerate}
Since $\kappa_2<\gb$, we may choose $g\in {}^\omega\omega$ such that 
\[g^*<g\quad\mbox{ and }\quad (\forall\alpha<\kappa_2)(g^\alpha<^*g).\] 
Pick $\bar{n}\in S$ which obeys $g$ and consider the function
$h_{\bar{n}}$. It follows from $(\boxdot)_1$ that $h_{\bar{n}}<^* f_i$ for
all $i<\omega$, so by the properties of an $(\aleph_0,\kappa_2)$--peculiar
cut there is $\alpha<\kappa_2$ such that $h_{\bar{n}}\leq^* f^\alpha$. Then,
for sufficiently large $i<\omega$ we have 
\begin{itemize}
\item $h_{\bar{n}}\restriction [n_i,n_{i+1})=f_i\restriction [n_i,n_{i+1})
    \leq f^\alpha\restriction [n_i,n_{i+1})$, and 
\item $n_i<g^\alpha(n_i)<g(n_i)<n_{i+1}$.
\end{itemize}
The latter implies that for some $\ell\in [n_i,n_{i+1})$ we have
  $f^\alpha(\ell)<f_i(\ell)$, contradicting the former. 
\end{proof}

\begin{theorem}
\label{y.10} 
Assume $\gp<\gt$. Then for some regular cardinal $\kappa$ there exists a
$(\kappa,\gp)$--peculiar cut in ${}^\omega\omega$ and $\aleph_1\leq
\kappa<\gp$.
\end{theorem}

\begin{proof}
Use Proposition \ref{y.2} and Lemma \ref{y.9} to choose
$\cB,\kappa,\bar{A},\pr$ and $\langle\bar{\eta}^\alpha:\alpha\le
\gp\rangle$ so that: 
\begin{enumerate}
\item[(i)] ${\mathcal B}=\{B_\alpha:\alpha<\gp\}$ exemplifies $\gp$,
\item[(ii)] $\bar{A}=\langle A_i:i<\kappa\rangle\subseteq
  [\omega]^{\aleph_0}$ is $\subseteq^*$--decreasing, $\kappa=\cf(\kappa) <
  \gp$ and conditions (a)+(b) of \ref{y.2} hold,
\item[(iii)]  $\pr:\gp\times\gp\longrightarrow\gp$ is a bijection satisfying
  $\pr(\alpha_1,\alpha_2)\ge\alpha_1,\alpha_2$,     
\item[(iv)] the sequence $\langle\bar{\eta}^\alpha:\alpha\le
\gp\rangle$ satisfies conditions (a)--(e) of \ref{y.9}.
\end{enumerate}
It is enough to find a suitable cut $\langle f_i:i<\kappa\rangle,\langle 
f^\alpha:\alpha<\gp\rangle\subseteq {}^{A^*}\omega$ for some infinite
$A^* \subseteq \omega$ (as by renaming, $A^*$ is $\omega$).  Let
\begin{enumerate}
\item[(v)] $A^*= \dom(\bar{\eta}^{\gp})$, 
\item[(vi)] for $i<\kappa$ we let $f_i:A^*\longrightarrow\omega$ be defined
  by  
\[f_i(n)=\min\big\{\ell:[\eta^{\gp}_n(n+\ell) =1\ \wedge n + \ell\notin
A_i]\mbox{ or }\dom(\eta^{\gp}_n)= [n,n + \ell)\big\},\]
\item[(vii)] for $\alpha<\gp$ we let $f^\alpha:A^*\longrightarrow\omega$
  be defined by 
\[f^\alpha(n) =\min\big\{\ell+1:[\eta^{\gp}_n(n+\ell)=1\ \wedge\ n+\ell\in
B_\alpha]\mbox{  or }\dom(\eta^{\gp}_n) = [n,n+\ell)\big\}.\] 
\end{enumerate}
Note that (by the choice of $f_i$, i.e., clause (vi)):
\begin{enumerate}
\item[(viii)]  $\bigcup\{[n,n + f_i(n)) \cap\set(\eta^{\gp}_n):n\in A^*\}
  \subseteq^* A_i$ for every $i < \kappa$. 
\end{enumerate}
Also,
\begin{enumerate}
\item[$(\circledast)^{\rm a}_1$] $f_j\leq^* f_i$ for $i<j<\kappa$.  
\end{enumerate}
[Why?  Let $i<j<\kappa$. Then $A_j \subseteq^* A_i$ and hence for some $n^*$
we have that $A_j\setminus n^* \subseteq A_i$. Therefore, for every $n\in A^*
\setminus n^*$, in the definition of $f_i,f_j$ in clause (vi), if $\ell$ can
serve as a candidate for $f_i(n)$ then it can serve for $f_j(n)$, so (as we
use the minimum there) $f_j(n)\le f_i(n)$.  Consequently $f_j \le^* f_i$.]

Now, we want to argue that we may find a subsequence of $\langle
f_i:i<\kappa\rangle$ which is $<^*$--decreasing. For this it is enough to
show that 
\begin{enumerate}
\item[$(\circledast)^{\rm b}_1$] for every $i<\kappa$, for some $j\in
  (i,\kappa)$ we have $f_j <^* f_i$. 
\end{enumerate}
So assume towards contradiction that for some $i(*)<\kappa$ we have
\[(\forall j)(i(*) < j < \kappa\ \Rightarrow\ \neg(f_j <^* f_{i(*)})).\]
For $j<\kappa$ put $B^*_j=:\{n\in A^*:f_j(n) \ge f_{i(*)}(n)\}$. Then $B^*_j
\in [A^*]^{\aleph_0}$ is $\subseteq^*$--decreasing, so there is a
pseudo-intersection $B^*$ of $\langle B^*_j:j <\kappa \rangle$ (so $B^* \in
[A^*]^{\aleph_0}$ and $(\forall j<\kappa)(B^*\subseteq^* B^*_j)$).  Now,
let $A'=\bigcup\{\set(\eta^{\gp}_n)\cap [n,n + f_{i(*)}(n)):n \in B^*\}$. 

\begin{enumerate}
\item[$(*)$] $A'$ is an infinite subset of $\omega$.
\end{enumerate}
[Why?  Recall that by \ref{y.9}(a) we have $\bar{\eta}^0 \in
{\mathbf{S}}_{\bar A}$ and hence $\set(\bar{\eta}^0)\subseteq^* A_{i(*)}$
and $(\forall n\in\dom(\bar{\eta}^0)(\set(\eta_n^0)\neq\emptyset)$ (see
Definition \ref{y.5}(2)). By \ref{y.9}(b) we know that for every large 
enough $n\in \dom(\bar{\eta}^{\gp})$ we have $n\in\dom(\bar{\eta}^0)$  and  
$\eta^0_n \trianglelefteq \eta^{\gp}_n$.  For every large enough $n\in
\dom(\bar{\eta}^0)$ we have $\set(\bar{\eta}^0) \setminus \{0,\ldots,n-1\}
\subseteq A_{i(*)}$,  and hence for every large enough $n\in
\dom(\bar{\eta}^{\gp})$ we have $\eta^0_n\trianglelefteq \eta^{\gp}_n$
and $\emptyset \neq\set(\eta^0_n) \subseteq A_{i(*)}$. Consequently, for
large enough $n\in B^*$, $[n,n+f_{i(*)}(n)) \cap\set(\eta^{\gp}_n) \ne
\emptyset$ and we are done.]  

\begin{enumerate}
\item[$(**)$] $A' \subseteq^* A_j$ for $j \in (i(*),\kappa)$ (and hence
for all $j < \kappa$).
\end{enumerate} 
[As $f_j \restriction B^* =^* f_{i(*)} \restriction B^*$ for $j \in (i(*),\kappa)$.] 

\begin{enumerate}
\item[$(***)$] $A' \cap B_\alpha$ is infinite for $\alpha<\gp$.  
\end{enumerate}
[Why?  By clauses (c) + (a) of \ref{y.9}, for every large enough $n\in
\dom(\bar\eta^{\alpha +1})$ we have set$(\eta^{\alpha +1}_n) \cap  B_\alpha
\ne \emptyset$ and set$(\eta^{\alpha +1}_n) \subseteq A_{i(*)}$.] 

\noindent Properties $(*)$--$(***)$ contradict (b) of \ref{y.2}, finishing
the proof of $(\circledast)^{\rm b}_1$.
\medskip

Thus passing to a subsequence if necessary, we may assume that 
\begin{enumerate}
\item[$(\circledast)^{\rm c}_1$] the demand in $(\alpha)$ of \ref{cut} is
  satisfied, i.e., $f_j<^* f_i$ for $i<j<\kappa$.  
\end{enumerate}
\medskip

Now,
\begin{enumerate}
\item[$(\circledast)_2$] $(\forall i<\kappa)(\forall\alpha<\gp)(f^\alpha <^*
  f_i)$.  
\end{enumerate}
[Why? Let $i<\kappa$, $\alpha< \gp$. For large enough $n\in A^*$ we have
that $\set(\eta^{\alpha+1}_n)\subseteq A_i$ and $\set(\eta^{\alpha+1}_n)\cap
B_\alpha\neq \emptyset$ and $\eta^{\alpha+1}_n\trianglelefteq
\eta^{\gp}_n$. Then for those $n$ we have $f^\alpha(n)\leq f_i(n)$. Now,
remembering $(\circledast)_1$, we may conclude that actually $f^\alpha<^*
f_i$.]  
\medskip

\begin{enumerate}
\item[$(\circledast)_3^{\rm a}$] The set (of functions) $\{f_i:i<\kappa\}
  \cup \{f^\alpha:\alpha<\gp\}$ is linearly ordered by $\le^*$,
\item[$(\circledast)_3^{\rm b}$] in fact, if $f',f''$ are in the family then
  either $f' =^* f''$ or $f'<^* f''$ or $f''<^* f'$. 
\end{enumerate}
[Why?  By $(\circledast)_1$, $(\circledast)_2$ and clauses (d) + (e) of
\ref{y.9}.]  
\medskip

Choose inductively a sequence $\bar{\alpha}=\langle\alpha(\varepsilon):
\varepsilon<\varepsilon^* \rangle\subseteq\gp$ such that: 
\begin{itemize}
\item $\alpha(\varepsilon)$ is the minimal $\alpha\in\gp\setminus\{ \alpha
  (\zeta):\zeta < \varepsilon\}$ satisfying $(\forall\zeta<\varepsilon)
  (f^{\alpha(\zeta)} <^* f^\alpha)$, and  
\item we cannot choose $\alpha(\varepsilon^*)$.
\end{itemize}
We ignore (till $(\circledast_7)$) the question of the value of $\varepsilon^*$.  Now,
\medskip

\begin{enumerate}
\item[$(\circledast)_4$]  $\langle f_i:i < \kappa \rangle$, $\langle
  f^{\alpha(\varepsilon)}: \varepsilon < \varepsilon^*\rangle$
satisfy clauses $(\alpha)$--$(\gamma)$ of \ref{cut}.
\end{enumerate}
[Why?  By $(\circledast)_1$--$(\circledast)_3$ and the choice of
$\alpha(\varepsilon)$'s above.]
\medskip

\begin{enumerate}
\item[$(\circledast)_5$] $\langle f_i:i<\kappa \rangle$, $\langle
  f^{\alpha(\varepsilon)}:\varepsilon<\varepsilon^*\rangle$ satisfy
  clause $(\delta)$ of \ref{cut}.
\end{enumerate}
[Why?  Assume towards contradiction that $f:A^*\longrightarrow \omega$ and
\[\big(\forall i<\kappa\big)\big(f\le^* f_i\big)\ \mbox{ but }\ \big(\forall
\varepsilon<\varepsilon^*\big)\big(\neg(f \le^* f^{\alpha(\varepsilon)}) 
\big).\]  
Clearly, without loss of generality, we may assume that $[n,n+f(n))
\subseteq\dom(\eta^{\gp}_n)$ for $n\in A^*$. Let $A' =
\bigcup\big\{[n,n+f(n)) \cap\set(\eta^{\gp}_n):n \in A^*\big\}$.  Now for
every $i<\kappa$, $A' \subseteq^* A_i$ because $f \le^* f_i$ and by the
definition of $f_i$. Also, for every $\alpha<\gp$, the intersection $A' \cap
B_\alpha$ is infinite.  Why? It follows from the choice of the sequence
$\bar{\alpha}$ that for some $\varepsilon<\varepsilon^*$ we have $\neg
(f^{\alpha(\varepsilon)}<^* f^\alpha)$ , and thus $f^\alpha\leq^*
f^{\alpha(\varepsilon)}$ (remember $(\circledast)_3$). Hence, if $n \in A^*$
is large enough, then $f^\alpha(n) \le f^{\alpha(\varepsilon)}(n)$ and for
infinitely many $n \in A^*$ we have $f^\alpha(n) \le f^{\alpha(\varepsilon)}
(n) < f(n) \le f_0(n) \le |\dom(\eta^{\gp}_n)|$. For every such $n$ we have
$n+f^\alpha(n)-1\in A' \cap B_\alpha$. Together,  $A'$ contradicts clause
(ii) of the choice of $\langle A_i:i<\kappa\rangle$, $\langle B_\alpha:
\alpha<\gp\rangle$, specifically the property stated in \ref{y.2}(b).]  
\medskip 

\begin{enumerate}
\item[$(\circledast)_6$]  $\langle f_i:i < \kappa \rangle,\langle 
f^{\alpha(\varepsilon)}:\varepsilon < \varepsilon^*\rangle$ satisfy clause
$(\varepsilon)$ of \ref{cut}.  
\end{enumerate}
[Why?  Assume towards contradiction that $f:A^*\longrightarrow\omega$, and
\[\big(\forall\varepsilon<\varepsilon^*\big)\big(f^{\alpha(\varepsilon)}
\le^* f\big)\ \mbox{ but }\ \big(\forall i<\kappa\big)\big(\neg(f_i \le^* 
f)\big).\]   
It follows from $(\circledast)_1$ (and the assumption above) that we may
choose an infinite set $A^{**} \subseteq A^*$ such that $\big(\forall
i<\kappa\big)\big((f \restriction A^{**}) <^* (f_i\restriction
A^{**})\big)$.  Let  
\[A'=\bigcup\{[n,n + f(n))\cap\set(\eta^{\gp}_n):n \in A^{**}\}\subseteq 
\omega.\] 
Since $(f \restriction A^{**}) <^* (f_i\restriction A^{**})$ we easily see
that $A' \subseteq^* A_i$ for all $i<\kappa$ (remember (viii)). Like in the
justification for $(\circledast)_5$ above, if $\alpha<\gp$ then for some
$\varepsilon<\varepsilon^*$ we have $f^\alpha \le^*f^{\alpha(\varepsilon)}$
and we may conclude from our assumption towards contradiction that $f^\alpha
\le^* f$ for all $\alpha<\gp$. Like in $(\circledast)_5$ we conclude that
for every $\alpha<\gp$ the intersection $A' \cap B_\alpha$ is infinite,
contradicting the choice of $\langle A_i:i<\kappa\rangle$, $\langle
B_\alpha:\alpha < \gp \rangle$.]

\begin{enumerate}
\item[$(\circledast)_7$] $\varepsilon^*=\gp$.
\end{enumerate}
[Why?  The sequence $\langle\alpha(\varepsilon):\varepsilon<\gp\rangle$ is
an increasing sequence of ordinals $<\gp$,  hence $\varepsilon^*\le \gp$. If
$\varepsilon^*<\gp$, then by the Bell theorem we get contradiction to
$(\circledast)_4$--$(\circledast)_6$ above; cf.~Proposition \ref{y.17} below.]   
\medskip

So $\langle f_i:i < \kappa \rangle,\langle f^{\alpha(\varepsilon)}:
\varepsilon<\gp\rangle$ are as required: clauses $(\alpha)$--$(\gamma)$ of
\ref{cut} hold by $(\circledast)_4$, clause $(\delta)$ by $\circledast_5$,
clause $(\varepsilon)$ by $(\circledast)_6$. Finally, since $\gt\leq \gb$,
we may use Proposition \ref{getb} to conclude that (under our assumption
$\gp<\gt$) there is no $(\aleph_0,\gp)$--peculiar cut and hence
$\kappa\geq\aleph_1$.
\end{proof}

\begin{remark}
  The existence of $(\kappa,\gp)$--peculiar cuts for $\kappa<\gp$ is
  independent from ``ZFC+$\gp=\gt$''. We will address this issue in
  \cite{Sh:F769}.
\end{remark}

\section{Peculiar cuts and {\bf MA}}

\begin{proposition}
\label{y.17} 
Assume that $\kappa_1\le\kappa_2$ are infinite regular cardinals and there
exists a $(\kappa_1,\kappa_2)$--peculiar cut in ${}^\omega\omega$. Then for
some $\sigma$--centered forcing notion $\bbQ$ of cardinality $\kappa_1$ and
a sequence $\langle\cI_\alpha:\alpha <\kappa_2\rangle$ of open dense subsets
of  $\bbQ$, there is no directed $G\subseteq\bbQ$ such that $(\forall\alpha<
\kappa_2)(G\cap\cI_\alpha\ne\emptyset)$. Hence ${\bf MA}_{\kappa_2}(
\sigma\mbox{--centered})$ fails and thus $\gp\leq\kappa_2$. 
\end{proposition}

\begin{proof}
Let $\big(\langle f_i:i<\kappa_1\rangle,\langle f^\alpha:\alpha<\kappa_2
\rangle\big)$ be a $(\kappa_1,\kappa_2)$--peculiar cut in
${}^\omega\omega$. Define a forcing notion $\bbQ$ as follows. 

\noindent {\bf A condition in $\bbQ$} is a pair $p= (\rho,u)$ such that
$\rho\in {}^{\omega >}\omega$ and $u\subseteq\kappa_1$ is finite.

\noindent {\bf The order $\leq_{\bbQ}=\leq$ of $\bbQ$} is given by
$(\rho_1,u_1)\leq (\rho_2,u_2)$ if and only if (both are in $\bbQ$ and) 
\begin{enumerate}
\item[(a)] $\rho_1\trianglelefteq \rho_2$,
\item[(b)] $u_1 \subseteq u_2$, 
\item[(c)]  if $n\in [\lh(\rho_1),\lh(\rho_2))$ and $i\in u_1$,  then $f_i(n)\geq \rho_2(n)$.
\end{enumerate}

Plainly, $\bbQ$ is a forcing notion of cardinality $\kappa_1$. It is
$\sigma$--centered as for each $\rho\in {}^{\omega}\omega$, the set
$\{(\eta,u)\in\bbQ: \eta=\rho\}$ is directed.  

For $j<\kappa_1$ let $\cI_j=\{(\rho,u) \in\bbQ:j\in u\}$, and for $\alpha=
\omega \beta+n<\kappa_2$ let 
\[\cI^\alpha=\big\{(\rho,u)\in\bbQ:\big(\exists m<\lh(\rho)\big)\big(m\ge n\
\wedge \rho(m)> f^\beta(m)\big)\big\}.\]
Clearly $\cI_j,\cI^\alpha$ are dense open subsets of $\bbQ$.  Suppose
towards contradiction that there is a directed $G\subseteq\bbQ$ intersecting 
all $\cI^\alpha,\cI_j$ for $j<\kappa_1$, $\alpha<\kappa_2$.  Put $g= 
\bigcup\{\rho:(\exists u)((\rho,u) \in G)\}$. Then 
\begin{itemize}
\item $g$ is a function, its domain is $\omega$ (as $G\cap \cI^n\ne
  \emptyset$ for $n<\omega$), and  
\item $g\le^* f_i$ (as $G \cap \cI_i\ne \emptyset$), and  
\item $\{n<\omega:f^\alpha(n)<g(n)\}$ is infinite (as $G\cap \cI^{\omega
    \alpha+n}\ne\emptyset$ for every $n$). 
\end{itemize}
The properties of the function $g$ clearly contradict our assumptions on 
$\langle f_i:i<\kappa_1\rangle$, $\langle f^\alpha:\alpha<\kappa_2\rangle$. 
\end{proof}

\begin{corollary}
If there exists an $(\aleph_0,\kappa_2)$--peculiar cut, then ${\rm
  cov}(\cM)\leq \kappa_2$.
\end{corollary}

\begin{theorem}
\label{y.27}  
Let $\cf(\kappa_2)=\kappa_2>\aleph_1$. Assume ${\bf MA}_{\aleph_1}$
holds. Then there is no $(\aleph_1,\kappa_2)$--peculiar cut in
${}^\omega\omega$. 
\end{theorem}

\begin{proof}
Suppose towards contradiction that $\cf(\kappa_2)=\kappa_2>\aleph_1$,
$\big(\langle f_i:i<\omega_1\rangle,\langle f^\alpha:\alpha<\kappa_2\rangle
\big)$ is an $(\aleph_1,\kappa_2)$--peculiar cut and ${\bf MA}_{\aleph_1}$
holds true. We define a forcing notion $\bbQ$ as follows.   

\noindent {\bf A condition in $\bbQ$} is a pair $p= (u,\bar{\rho})=
(u^p,\bar{\rho}^p)$ such that 
\begin{enumerate}
\item[(a)] $u\subseteq\omega_1$ is finite, $\bar{\rho}=\langle\rho_i:i \in u
  \rangle=\langle\rho_i^p:i \in u\rangle$,  
\item[(b)] for some $n=n^p$, for all $i \in u$ we have $\rho_i\in
  {}^n\omega$,  
\item[(c)]  for each $i\in u$ and $m<n^p$ we have $\rho_i(m)\leq f_i(m)$, 
\item[(d)] if $i_0=\max(u)$ and $m\geq n^p$, then $f_{i_0}(m)>2\cdot
  |u^p|+885$. 
\item[(e)] $\langle f_i \restriction [n^p,\omega):i\in u\rangle$ is
  $<$--decreasing.
\end{enumerate}

\noindent {\bf The order $\leq_{\bbQ}=\leq$ of $\bbQ$} is given by
$p\leq q$ if and only if ($p,q\in\bbQ$ and) 

\begin{enumerate}
\item[(f)] $u^p \subseteq u^q$, 
\item[(g)] $\rho^p_i \trianglelefteq \rho^q_i$ for every $i \in u^p$, 	
\item[(h)]  if $i<j$ are from $u^p$, then $\rho^q_i\restriction [n^p,n^q)<
  \rho^q_j \restriction [n^p,n^q)$, 
\item[(i)] if $i<j$, $i\in u^q\setminus u^p$ and $j\in u^p$, then for some
  $m\in [n^p,n^q)$ we have $f_j(m)<\rho^q_i(m)$.
\end{enumerate}

\begin{claim}
\label{cl1}
$\bbQ$ is a ccc forcing notion of size $\aleph_1$.
\end{claim}

\begin{proof}[Proof of the Claim]
Plainly, the relation $\leq_\bbQ$ is transitive and $|\bbQ|=\aleph_1$. Let 
us argue that the forcing notion $\bbQ$ satisfies the ccc.

Let $p_\varepsilon\in\bbQ$ for $\varepsilon<\omega_1$. Without loss of
generality $\langle p_\varepsilon:\varepsilon<\omega_1\rangle$ is without
repetition.  Applying the $\Delta$--Lemma we can find an unbounded set
$\cU\subseteq\omega_1$ and $m(*)<n(*)<\omega$ and $n'<\omega$ such that for
each $\varepsilon\in\cU$ we have
\begin{enumerate}
\item[(i)]  $|u^{p_\varepsilon}|=n(*)$ and $n^{p_\varepsilon}=n'$; let
  $u^{p_\varepsilon}=\{\alpha_{\varepsilon,\ell}:\ell<n(*)\}$ and
  $\alpha_{\varepsilon,\ell}$ increases with $\ell$, and 
\item[(ii)] $\alpha_{\varepsilon,\ell} =\alpha_\ell$ for $\ell<m(*)$ and
  $\rho_{\varepsilon,\ell} = \rho^*_\ell$ for $\ell <  n(*)$, and 
\item[(iii)] if $\varepsilon<\zeta$ are from $\cU$ and $k,\ell \in
  [m(*),n(*))$, then $\alpha_{\varepsilon,\ell}<\alpha_{\zeta,k}$.
\end{enumerate}
Let $\varepsilon<\zeta$ be elements of $\cU$ such that
$[\varepsilon,\zeta)\cap\cU$ is infinite. Pick $k^*>n'$ such that for each
$k\geq k^*$ we have 
\begin{itemize}
\item the sequence $\langle f_\alpha(k):\alpha\in\{\alpha_{\varepsilon,
    \ell}:\ell<n(*)\}\cup\{\alpha_{\zeta,\ell}:\ell<n(*)\}\rangle$ is strictly decreasing, 
\item $f_{\alpha_{\zeta,n(*)-1}}(k)>885\cdot (n(*)+1)$, and 
\item $f_{\alpha_{\zeta,m(*)}}(k)+n(*)+885< f_{\alpha_{\varepsilon,n(*)-1}}
  (k)$.  
\end{itemize}
(The choice is possible because $\langle f_i:i<\omega_1\rangle$ is
$<^*$--decreasing and by the selection of $\varepsilon,\zeta$ we also have
$\lim\limits_{k\to\infty} \big(f_{\alpha_{\varepsilon,n(*)-1}}(k)- 
f_{\alpha_{\zeta,m(*)}} (k)\big)=\infty$.)

Now  define $q = (u^q,\bar{\rho}^q)$ as follows:  
\begin{itemize}
\item $u^q = u^{p_\varepsilon} \cup u^{p_\zeta}$, $n^q=k^*+1$, 
\item if $n<n'$, $i\in u^{p_\varepsilon}$, then
  $\rho^q_i(n)=\rho^{p_\varepsilon}_i(n)$, 
\item if $n<n'$, $i\in u^{p_\zeta}$, then $\rho^q_i(n)=\rho^{p_\zeta}_i(n)$,  
\item if $i=\alpha_{\varepsilon,\ell}$, $\ell<n(*)$, $n\in [n',k^*)$, then
  $\rho^q_i(n)=\ell$, and if $j=\alpha_{\zeta,\ell}$, $m(*)\leq\ell<n(*)$,
  then $\rho^q_j(n)=n(*)+\ell+1$, 
\item if $j=\alpha_{\zeta,\ell}$, $\ell<n(*)$, then $\rho^q_j(k^*)=\ell$,
  and if $i=\alpha_{\varepsilon,\ell}$, $m(*)\leq \ell<n(*)$, then
  $\rho^q_i(k^*)=f_{\alpha_{\zeta,m(*)}}(k^*)+\ell+1$. 
\end{itemize}
It is well defined (as $\rho^{p_\varepsilon}_{\alpha_{\varepsilon,\ell}} =
\rho^{p_\zeta}_{\alpha_{\zeta,\ell}}$ for $\ell < m(*))$. Also $q \in
\bbQ$. Lastly, one easily verifies that $p_\varepsilon\le_{\bbQ} q$ and
$p_\zeta\le_{\bbQ} q$, so indeed $\bbQ$ satisfies the ccc.
\end{proof}

For $i < \omega_1$ and $n < \omega$ let 
\[\cI_{i,n}=\big\{p\in\bbQ:\big[u^p\nsubseteq i \mbox{ or for no $q\in\bbQ$ we
  have }p \le_{\bbQ} q \wedge u^q \nsubseteq i\big]\mbox{ and } n^p \ge
n\big\}.\]   
Plainly, 
\begin{enumerate}
\item[$(\alpha)$] the sets $\cI_{i,n}$ are open dense in $\bbQ$. 
\end{enumerate}
Also,
\begin{enumerate}
\item[$(\beta)$] for each $i<\omega_1$ there is $p^*_i \in \bbQ$ such that 
  $u^{p_i}=\{i\}$. 
\end{enumerate}
It follows from \ref{cl1} that 
\begin{enumerate}
\item[$(\delta)$] for some $i(*)$, $p^*_{i(*)} \Vdash_{\bbQ}$
  ``$\{j<\omega_1: p^*_j\in \name{G}\}$ is unbounded in $\omega_1$ ''. 
\end{enumerate}
Note also that
\begin{enumerate}
\item[$(\gamma)$] if $p$ is compatible with $p^*_{i(*)}$ and $p\in
  \cI_{i,n}$ then $u_p \nsubseteq i$.  
\end{enumerate}

Since we have assumed ${\bf MA}_{\aleph_1}$ and $\bbQ$ satisfies the ccc (by
\ref{cl1}), we may find a directed set $G\subseteq\bbQ$ such that
$p^*_{i(*)} \in G$ and $\cI_{i,n} \cap G \ne\emptyset$ for all $n<\omega$
and  $i<\omega_1$. Note that then the set $\cU := \bigcup\{u^p:p \in G\}$ is
unbounded in $\omega_1$.

For $i \in \cU$ let $g_i=\bigcup\{\rho^p_i:p \in G\}$. Clearly each $g_i\in
{}^\omega \omega$ (as $G$ is directed, $\cI_{i,n}\cap G \ne \emptyset$ for
$i < \omega_1$, $n<\omega$). Also $g_i \le f_i$ by clause (c) of the
definition of $\bbQ$, and $\langle g_i:i \in \cU\rangle$ is
$<^*$--increasing by clause (h) of the definition of $\le_{\bbQ}$.  Hence 
for each $i < j$ from $\cU$ we have $g_i <^* g_j\le^* f_j <^* f_i$. Thus by
property \ref{cut}$(\delta)$ of a peculiar cut, for every $i\in\cU$ there is
$\gamma(i)<\kappa_2$ such that $g_i <^* f^{\gamma(i)}$.  Let
$\gamma(*)=\sup\{\gamma(i):i\in\cU\}$. Then $\gamma(*) < \kappa_2$ (as   
$\kappa_2=\cf(\kappa_2)>\aleph_1$). Now, for each $i\in\cU$ we have $g_i<^*
f^{\gamma(*)}<^*f_i$, and thus for $i\in\cU$ we may pick $n_i<\omega$ such that   
\[n\in [n_i,\omega)\ \Rightarrow\ g_i(n)< f^{\gamma(*)}(n)<f_i(n).\]   
For some $n^*$ the set $\cU_* = \{i\in \cU:n_i= n^*\}$ is unbounded in
$\omega_1$.  Let $j\in \cU_*$ be such that $\cU_* \cap j$ is infinite. Pick
$p\in G$ such that $j\in u^p$ and $n^p>n^*$ (remember $G\cap
\cI_{j,{n^*+1}}\neq \emptyset$ and $G$ is directed). Since $u^p$ is finite,
we may choose $i\in\cU_*\cap j\setminus u^p$, and then $q\in G$ such that 
$q\geq p$ and $i\in u^q$. If follows from clause (i) of the definition of
the order $\leq$ of $\bbQ$ that there is $n\in [n^p,n^q)$ such that
$f_j(n)<\rho^q_i(n)=g_i(n)$. Since $n>n^*=n_i=n_j$ we get 
$f_j(n)<g_i(n)<f^{\gamma(*)}(n)<f_j(n)$,  a contradiction.   
\end{proof}

\begin{remark}
The proof of \ref{y.27} actually used Hausdorff gaps on which much is
known (see, e.g., Abraham and Shelah \cite{AbSh:537}, \cite{AbSh:598}). More
precisely, the proof could be presented as a two-step argument:
\begin{enumerate}
\item from ${\bf MA}_{\aleph_1}$ one gets that every decreasing
  $\omega_1$--sequence is half-a-Hausdorff gap, and
\item if $\kappa_2=\cf(\kappa_2)>\aleph_1$, then the $\omega_1$--part of a
  peculiar $(\omega_1,\kappa_2)$--cut cannot be half-a-Hausdorff gap.
\end{enumerate}
\end{remark}

\begin{corollary}
\label{y.49} 
If ${\bf MA}_{\aleph_1}$ then $\gp=\aleph_2\ \Leftrightarrow \gt=\aleph_2$.  In
other words, 
\[\gm=\gp=\aleph_2\ \Rightarrow\ \gt=\aleph_2.\]  
\end{corollary}

\bibliographystyle{plain}
\bibliography{lista,listb,listx,listf,liste,listy}

\end{document}